\documentstyle{article}
\input{amssym.def}
\input{amssym.tex}

\newcommand{\forces}{\; \Vdash}

\newcommand{\Pforces}{\ \Vdash_P\ }
\newcommand{\Qforces}{\ \Vdash\ }										
\newcommand{\Kforces}{\ \Vdash_K\ }										
\newcommand{\Rforces}{\ \Vdash\ }										
\newcommand{\rest}{\restriction}
\newcommand{\ltree}{\lambda ^{< \omega}}
\newcommand{\powertaumu}{{ P}^\tau_\mu}

\newcommand{\A}{Aronszajn}
\newcommand{\ga} {\alpha}
\newcommand{\gb}{\beta}	
\newcommand{\gd} {\delta}
\newcommand{\gr} {\rho}
\newcommand{\gs} {\sigma}
\newcommand{\gl} {\lambda}
\newcommand{\go} {\omega}
\newcommand{\gz} {\zeta}
\newcommand{\gi} {\iota}
\newcommand{\gk} {\kappa}
\newcommand{\gt} {\tau}

\newcommand{\bfB}{{\bf B}}

\newcommand{\calU}{{\cal U}}

\newcommand{\alephomp}{\mbox{$\aleph_{\omega+1}$}}

\newcommand{\alephom}{\aleph_\omega}

\newcommand{\kpom}{\mbox{$\kappa^{+\omega}$}}
\newcommand{\kpompo}{\mbox{$\kappa^{+(\omega+1)}$}}

\newcommand{\jkpompo}{j(\kappa)^{+(\omega+1)}}

\newcommand{\la}{\langle}
\newcommand{\ra}{\rangle}

 \newcommand{\CS} {\mbox{$\cal S$}}

 \newcommand{\tildaT}{{\bf T}}
 
 \newcommand{\tildaR}{{\bf  R}}
  
\newcommand{\tildaCS}{\raisebox{-.2 cm} {$\stackrel{\CS}{\sim}$}}

 \newcommand{\cf}{\it cf}

 \newcommand{\seqlambda}{\la \gl_i \mid i < \go \ra}

 \newcommand{\Pre}{{\it Pre}}
\newcommand{\Coll}{{\it Coll}}
\newcommand{\dom}{{\it dom}}

\def\newtheorems{\newtheorem{theorem}{Theorem}[section]
                 \newtheorem{cor}[theorem]{Corollary}
                 
                 \newtheorem{lemma}[theorem]{Lemma}

                 \newtheorem{definition}{Definition}[section]}

\newtheorems

\begin{document}

\title{The tree property at successors of singular cardinals}
\author{Menachem Magidor and Saharon Shelah
\thanks{Partially sponsored by the Edmund Landau Center
for research in Mathematical Analysis, supported by the Minerva Foundation
(Germany), p.n. }\\
Institute of Mathematics\\
The Hebrew University, Jerusalem, Israel}
\maketitle
\begin{abstract} Assuming some large cardinals, a model of ZFC is obtained
in which \alephomp\ carries no \A\ trees. It is also shown that if $\gl$
is a singular limit of strongly compact cardinals, then $\gl^+$ carries no
\A\ trees.
\end{abstract}

\section{Introduction}
%$\frak P$
\label{s1}
The main results of this paper are (1) that the consistency of ``ZFC and
\alephomp\ carries no \A\ trees'' follows from the consistency of some large
cardinals (a huge cardinal with $\go$ supercompact cardinals
above it), and (2) that if a singular cardinal $\gl$ is
a limit of strongly compact cardinals, then there are no \A\ trees of height
$\gl^+$. The proof of (2) is in Section 3,
and the forcing constructions which prove (1) are given in Sections 4, 5,
and 6.
  The generalization to higher singular cardinals of both
(1) and (2) poses no
problem.

Section \ref{s2} contains
 some notations, definitions, and standard
constructions.

\section{Preliminaries}
\label{s2}
We use the convention concerning forcing by which $p < q$ means that $q$ is
more informative than $p$. A forcing poset for us here is a separative
partial order $P$, with a least informative point (denoted $\emptyset_P$) and
 with no maximal point. (A poset is separative if:
Whenever  $p \not
< q$ then some extension of $q$ in $P$ is incompatible with $p$.
There is a canonical way of producing a separative poset from
non-separative: Define
an equivalence relation $p_1 \sim p_2$
if
  ``any $x$ is compatible with $p_1$ iff it is compatible with $p_2$''.
  Then, on the equivalence classes,
 define $[p_1] \leq [p_2]$ if every $x$ compatible
  with $p_2$ is also compatible with $p_1$.)

$V^P$ denotes the class of all $P$--terms, but when an expression such as
``${\bf T} \in V^P$ is a tree of hight $\lambda^+$'' is used, we mean that
$\emptyset_P$ forces this statement.

An \A\ tree of height $\gl^+$ ($\gl$ a cardinal) is a tree with $\gl^+$
levels, each of cardinality $\leq \gl$, such that there is no branch of
size $\gl^+$. For the question that interests us, it does no harm to assume
that {\em all} the levels are of size $\gl$, and it is technically useful
to assume that the $\ga$th level is always the set $\gl \times \{\ga\}$.
So if $K$ is a forcing poset and $\tildaT \in V^K$ a name of a
$\gl^+$-tree, then the underlying universe of $\tildaT$ is in $V$, but,
of course, the
ordering $<_{\tildaT}$ is obtained by forcing.
\begin{definition}
\label{pretree}
Let $\tildaT$ be a $K$-name of a $\gl^+$-tree as above, where $\gl$ and
$\gl^+$ are cardinals both in $V$ and $V^K$, then the {\em
pre-tree} of $\tildaT$ is the system of relations $\la R_p \mid p \in K
\ra$ satisfying
\[ a\ R_p \ b \mbox{ iff } p \Kforces a <_{\tildaT} b .\]
So $R_p$ is a relation on the universe of $T$, namely on $\lambda \times
\lambda^+$.
\end{definition}
It turns out that the consistency proof for ``no \A\ trees on $\lambda^+$"
 relies on an
investigation of such systems of pre-trees, when $\| K \| < \gl$; the
crucial properties of which
are collected in the following definition.
The reader who prefers to do so, however, can always think in terms of
trees, pre-trees etc. without any harm.
\begin{definition}
\label{d12}
\begin{description}
\item[Systems:]
Suppose that $\tau \leq \lambda$ are cardinals,
$D \subseteq \gl^+$ is unbounded, and 	
 $T=\la T_\alpha \mid \alpha \in D \ra$ is a sequence such that
$T_\alpha \subseteq \tau \times \{\alpha \}$, for $\alpha \in D$.
Let $I$ be an index set of
cardinality $\leq \lambda$, and	  $R = \{ R_i \mid i \in I \}$
 a collection of binary relations that are subsets    of
  $\bigcup \{ T_{\ga} \times T_{\gb}
\mid \ga < \gb \mbox{ in}\ D\}$. Then the pair $\CS=( T, R)$
is called {\em a system} over $\gl^+$ (or a $\lambda^+$--system)
if the following hold:
\begin{enumerate}
 \item For every $\ga < \gb$ in $D$, there are $a \in T_\ga,\ b \in
T_{\gb}$ and $i \in I$ such that $\la a,b \ra \in R_i$.
\item For every $i \in I$, and $\ga < \gb < \gamma$ in $D$,
if $a \in T_\alpha,\; b\in T_\beta,\; c \in T_\gamma$ are such that
$\la a,c \ra \in R_i$ and $\la b,c \ra \in R_i$, then $\la a,b \ra \in
R_i$.
\end{enumerate}
The set $D$ is called the domain of the system and $I$  its index
set.  The cardinal $\tau$ is the {\em width}, and $\gl^+$ the {\em
height} of the system.

An example of a $\lambda^+$--system is any $\lambda^+$--tree; in this
example  $R$ consists of a single relation--the tree ordering. The
pre-tree as defined above is more illustrative; the number of relations is
the cardinality of the forcing poset.
 \item[Strong systems:] If condition (1) above is replaced by: ``For every
 $\alpha < \beta$ in $D$, for every $b \in T_\beta$,
there are $a \in T_\alpha$, and $i \in I$ such that $\la a , b \ra \in
R_i$'' (but (2) remains
unchanged) then the system $\CS$  is called {\em a strong system}.

A pre-tree relative to some forcing poset $K$ is
in fact a strong system.

\item[Subsystems:] If $\CS$ is a system then a
subsystem of $\CS$ is a system obtained by taking an
unbounded subset $D_0 \subseteq D$, an initial segment $\tau_0 \leq
 \tau$, a subset $I_0 \subseteq I$, and restricting $\CS$ to $\la
\tau_0 \times\{  \ga \} \mid \ga \in D_0 \ra$ and to the relations with
index in $I_0$.

By definition a subsystem is a system (though item {\it 2} is inherited
automatically from the system, item {\it 1} may not). A subsystem of a
strong system may no longer be strong.
\item[Narrow systems:] A system $\CS$	
 is said to be $(\gr, \gi)$--{\em narrow}
if $\tau < \gr$, and $\mid I \mid < \gi$.  That is, its width is less
than $\gr$ and its index set has size  $< \gi$.  A $\gl^+$--system is
said to be {\em narrow} iff it is $(\gl, \gl)$--narrow.

 \item[Branches:]
A ``branch'' of the system is a set $B$ such that for some $ i \in I$ for
all $a,b \in B\ \la a , b \ra \in R_i$ or $\la b , a \ra \in R_i$.

Thus, a branch of a pre-tree gives a set $B$ forced by a single condition
to be linearly ordered.
\item[Pre-systems]
Suppose that $Q$ is  a forcing poset and $\tildaCS=
 \la  \tildaT, \tildaR \ra$ is a name in $V^Q$ of a
system with domain $\gl^+$, width $\gl_0 \leq \gl$, and
index set some cardinal $\gt$. Then the {\em pre-system},
$\Pre_Q(\tildaCS)$ is defined as the following $\lambda^+$-system in $V$.
$\Pre_Q(\tildaCS)$ has as index set the product
$Q \times \tau$, its width remains $\lambda_0$,
 and its relations $R_{q,i}$ are defined
for $q \in Q$ and $i \in \tau$ by:
\[ \la a , b \ra \in R_{q,i}\ \mbox{iff}\	q \Qforces \la a , b \ra \in
\tildaR_i.\]

A pre-tree is an example of a pre-system. A pre-system is a system, and it
is strong if the system $\tildaCS$ is strong in $V^Q$ (that is, forced by
every condition to be strong).
If $P=Q \times K$ is a product and $\tildaCS$ is a $P$-name of a system,
then $\Pre_P(\tildaCS) \sim \Pre_Q(\Pre_K(\tildaCS))$.
\end{description}
\end{definition}

\subsection{On collapses}
For a regular cardinal $\kappa$ and an ordinal $\lambda > \kappa$,
$\Coll(\kappa, \lambda)$ is the poset that collapses $\lambda$ to
$\kappa$, using functions from ordinals below $\kappa$ and into
$\lambda$.  The poset $\Coll(\kappa, <\lambda)$ is the product with
support of cardinality $<\kappa$ of all collapses of ordinals between
$\kappa$ and $\lambda$ (The ``Levy'' collapse, see Jech, Section 20).
Now if $L=\la \lambda_i \mid i < \omega \ra$ is an increasing sequence
of cardinals, then $C=\Coll(L)$ is the full support iteration of the
collapsing posets $\Coll(\lambda_i, < \lambda_{i+1})$.  In more details,
define by induction on $n < \omega$ posets $P_n$ as follows.
$P_1=\Coll(\lambda_0,< \lambda_1)$, and
\[ P_{n+1}=P_n \ast \Coll(\lambda_{n-1}, < \lambda_n)^{V^{P_n}}.\]
Then $\Coll(L)$ is the full support iteration of the posets $P_n$.

A {\em projection} from a poset $P$ into $Q$ is an order preserving map
$\Pi : P \longrightarrow Q$ such that $\Pi(\emptyset_P)=\emptyset_Q$, and
 if $\Pi(p)=q$ and $q' > q$ in $Q$
then for some $p' > p$, $\Pi(p') \geq q'$.  (Some authors use a different
definition!)

If $\Pi : P \longrightarrow Q$ is a projection and $G \subset Q$ is a
$V$--generic filter, then $P/G$ is the separative poset defined by
taking $\Pi^{-1} G$ and  turning it into a separative poset.
  Then $P$ is isomorphic to a dense subset of the iteration $Q \ast
  (P/G)$: the isomorphism is the map taking $p \in P$ to $( \Pi(p),[p])$.

\begin{lemma} Suppose that $Q$ is a $\lambda$--complete poset (any
 increasing sequence of length $< \lambda$ has a {\em least} upper
bound).  Let $\mu$ be the cardinality of $Q$.  Then there is a
projection $\Pi : \Coll(\lambda, \mu) \longrightarrow Q$ such that
whenever $G \subset Q$ is $V$ generic, then the quotient poset
$\Coll(\lambda,\mu)/G$ is $\lambda$--complete.  In fact, the projection is
$\lambda$--continuous: If $f=\bigcup_{i< \lambda_0} f_i$ is the supremum
of an increasing sequence of lenght $\lambda_0 < \lambda$ of conditions
in $\Coll(\lambda,\mu)$ then $\Pi(f)$ is the supremum of $\{ \Pi(f_i)
\mid i < \lambda_0 \}$.
\end{lemma}
{\bf Proof}.  Let $Q=\{ q_i \mid i < \mu \}$ be an enumeration of $Q$.
Any condition in $\Coll(\lambda,\mu)$ is a function $f : \alpha
\longrightarrow \mu$, where $\alpha < \lambda$, and we define $\Pi(f)$ as
follows.  Define by induction an increasing sequence $\la a(\xi) \mid
\xi \leq \alpha \ra$ by requiring that (1) $a(0)$ is the minimum of $Q$
(2) at limit
stages $\delta \leq \alpha$, $a(\delta)$ is the least upper bound of
$\la a(\xi) \mid \xi < \delta \ra$, and (3) if $q_{f(\xi)}$ extends
$a(\xi)$ in $Q$, then $a(\xi+1)=q_{f(\xi)}$, and otherwise
$a(\xi+1)=a(\xi)$.  Finally, $\Pi(f)=a(\alpha)$.  It is easy to see that
$\Pi$ is a  $\lambda$-continuous projection.

Remark first that if $G
\subset Q$ is $V$-generic, then $\Pi^{-1}G$ is already separative. To prove
the $\lambda$--closure of the quotient,
suppose that $q \Qforces \la \tau_i \mid i< \lambda_0< \lambda \ra $
 is an increasing $Q$-sequence in  $\Coll(\lambda, \mu)/G= \Pi^{-1}G$.
We will find an extension of $q$ that forces a least upper bound to this
sequence.  Define by induction on $i<\lambda_0$ an increasing sequence
$\la q_i \mid i< \lambda_0 \ra$ beginning with $q_0=q$, such that for
every $i$, for
some $f_i \in \Coll(\lambda,\mu)$, $q_{i+1} \Qforces \tau_i = f_i$.
Let $q'$ be an upper bound in $Q$ to this sequence. Now
 $f_i \subset f_j$ for $i<j$, and $f=\bigcup_{i<
\lambda_0} f_i$ is a condition.  Since $Q$ is separative,
 $\Pi(f_i) \leq_Q q_{i+1}$ follows
from the fact that $q_{i+1}\forces \Pi(f_i)\in G$.  The continuity of
$\Pi$ implies that $\Pi(f)$ is the least upper bound of all the
conditions $\Pi(f_i)$, and hence $\Pi(f) \leq q'$.  That is $q'$ forces
that $f$ is in  $\Coll(\lambda, \mu)/G$.

A similar lemma holds for $\Coll(\lambda,<\kappa)$ if the cardinality of
$Q$ is less than $\kappa$ (and is $\lambda$--complete), but we need to apply
such a lemma in a slightly more complex situation.  Suppose that:
\begin{enumerate}
\item $\lambda < \lambda_1 \leq \mu < \kappa$ are regular cardinals.
$Q=\Coll(\lambda, < \lambda_1)$, and
the projection $\Pi_1: \Coll(\lambda, <\kappa) \longrightarrow Q$
 is the obvious restriction projection.
\item $P=Q \ast	R$ is a two stage iteration where
$R$ is a name in $V^{Q}$ such that
$R$ is forced to be $\lambda$--closed (by every condition).
 The projection of $P$ on $Q$ is denoted $\Gamma$ (so $\Gamma
(q,\tau)=q$).
% \item	The
%``term space'' for $R$ is the collection of all names $\tau \in
%V^{Q}$ that are forced by every condition to be in $R$, and where
%two names are identified if every condition forces them to be
%equal.	
%A partial ordering of the term space is defined when every condition forces
%$\tau_2$ to extend $\tau_1$.

 Suppose that the cardinality of  $P$ is $\mu$.
\end{enumerate}
\begin{lemma}
\label{l22}
Under the conditions set above on $P$, $Q$, and $R$, there is a
projection $\Pi : \Coll(\lambda, < \kappa) \longrightarrow P$ such
that $\Gamma \circ \Pi = \Pi_1$, and such that whenever $G \subseteq P$ is
$V$-generic, then the quotient poset $\Coll(\lambda, < \kappa)/G$ is
$\lambda$-closed.
\end{lemma}
{\bf Proof.} Set an enumeration $\{ \tau_i \mid i \in \mu \}$ of all the
terms
in $V^{Q}$ that are forced by every condition to be in $R$, and where
two names are identified if every condition forces them to be
equal.	
 Given any condition $q\in \Coll(\lambda, <
\kappa)$, let $q_1=\Pi_1(q)$, and $f=q\restriction \{ \mu \}$ be the
component of $q$ that collapses $\mu$.  Then $\Pi(q)=(q_1,\tau) \in Q \ast
R$ is
defined by the following procedure.  Suppose that the domain of $f$ is
$\alpha < \mu$, and define an increasing sequence $\la \eta_i \mid i \leq
\alpha \ra$ of terms, by induction on $i$ as follows:
\begin{enumerate}
\item $\eta_0$ is an assumed empty condition in $R$ (least informative).
$\eta_{i+1}$ is $\tau_{f(i)}$ if every condition in $Q$ forces that
$\tau_{f(i)}$ extends $\eta_i$, and $\eta_{i+1}=\eta_i$ otherwise.
\item If $\delta \leq \alpha$ is a limit ordinal and all the terms
$\eta_i$ for $i < \delta$ have been defined such that for $i<j<\delta$,
$\eta_i <_R \eta_j$ is forced by every condition in $Q$, then
$\eta_\delta$ is defined as (the name of) the least upper bound of this
increasing sequence.
\end{enumerate}
Finally, the projection is defined by setting $\tau=\eta_\alpha$.  We
leave it to the reader to verify that $\Pi$ is indeed a projection as
required, and in particular that if $G$ is generic over $P$ then the
quotient $\Coll(\lambda,<\kappa)/G$ is $\lambda$--closed.

\subsection{On embeddings and ultrapowers}
The dual characterization of supercompact cardinals is probably known to the
reader: If $\mu > \kappa$ then $\kappa$ is $\mu$-supercompact if $\kappa$
is the critical point (first fixed point)
 of an elementary embedding $j : V \longrightarrow M$
of the universe $V$ into a transitive inner model $M$ such that $\mu <
j(\kappa)$ and $M^\mu \subset M$. An alternative, more tangible, definition
is that $P_\kappa(\mu)=\{ X \mid X \subseteq \mu,\ | X | < \kappa \}$
carries a non-principal,
 fine, $\kappa$-complete, normal ultrafilter. (See Jech [] Chapter
6, or Solovay xxxxx).

We say that $\kappa$ is {\em huge} if $\kappa$ is the critical point of an
elementary embedding $j: V \longrightarrow M$ into a transitive
substructure $M$ such that $M^{j(\kappa)} \subseteq M$. If the stronger
demand, $M^{j(j(\kappa))} \subseteq M$, holds then we say that $\kappa$ is
2-huge. For our consistency result we need a cardinal that is slightly
stronger than huge, but not quite 2-huge. Its definition is given by
the following lemma
 on the
equivalence between two characterizations, which we quote without proof
(the proof is quite standard; see for example SRK).

If $\mu > \tau$ are cardinals, define $\powertaumu=\{ X
\subseteq \mu \mid \mbox{\it order-type of} \ X \ \mbox{\it is}\ \tau
\}$.  If ${\cal U} \subset {\cal P}(\powertaumu)$ is  a non-principal
ultrafilter, then ${\cal U}$ is said to be:
\begin{enumerate}
\item {\em $\kappa$--complete} if the intersection of fewer than
$\kappa $ sets in $\calU$ is again in $\calU$.
\item {\em normal} if any choice function $f$ is constant on a set in
$\calU$.  ($f$ is a choice function if $f(A) \in A$ for every $A \in
\powertaumu$.  $f$ is constant on $X \in \calU$ if for some $\gamma
\in \mu$, $f(A)=\gamma$ for all $A \in X$.)
\item {\em fine} if $\forall \alpha < \mu\ \{ A \mid \alpha \in A \}
\in \calU$.
\item {\em $\kappa$--small below} $\tau$ if $\{ A \in \powertaumu
\mid \  |A\cap \tau |< \kappa \} \in \calU$.
\end{enumerate}
\begin{lemma}
\begin{enumerate}
\item If $j:V \longrightarrow M$ is an elementary embedding into a
transitive structure $M$, with critical point $\kappa$, and $\tau$ is such
that $\kappa < \tau < j(\kappa) < j(\tau)=\mu$ and
 $M^\mu \subseteq M$, if $\calU \subset
{\cal P}(\powertaumu)$ is defined by
\[ A \in \calU\ \mbox{ iff}\ (j'' \tau) \in j(A),\]
then $\calU$ is a non-principal, $\kappa$--complete, normal, fine, and
$\kappa$--small below $\tau$ ultrafilter on $\powertaumu$.
\item If, on the other hand, $\kappa < \tau < \mu$ and $\calU$ are
such that $\calU \subset {\cal P}(\powertaumu)$ is an ultrafilter that satisfy
the five properties above, and if the ultrapower
$V^{\powertaumu}/\calU$ is created, then it is well--founded and the
ensuing elementary embedding, $i : V \longrightarrow N$, is such that
$\kappa$ is the critical point of $i$, $\tau < i(\kappa)$,
$\mu=i(\tau)$ and $N^{i(\tau)} \subset N$.
\end{enumerate}
\end{lemma}

%Let us give a name to the cardinal defined by this lemma, and say (only
%here, only here) that $\kappa$ is $\tau$-huge.2.3 if it satisfies
%the properties described in Lemma 2.3.
\subsection{The exact assumptions}
\label{s22}
For the consistency proof of ``no \A\ trees on $\alephomp$'' we need the
following:
\begin{quote}
 A cardinal $\kappa$ and an increasing sequence of cardinals $L=\la
\lambda_i \mid i < \omega \ra$ with $\lambda_0> \kappa$ such that:
\begin{description}
\item[$A_1$:] Let $\lambda= \sup \{ \lambda_i \mid i < \omega \}$ and $\mu =
\lambda^+$, then each $\lambda_i$, for $i > 0$, is $\mu$-supercompact.
\item[$A_2$:] If $P=\Coll(L)$, then in $V_1=V^P$ $\kappa$ is $\tau$-huge
for $\tau=\kpompo$ and an embedding $j : V_1 \longrightarrow M_1$ (as
in Lemma 2.3) such that  $j(\tau) = \mu$.
\end{description}
\end{quote}
This requirement of a ``potentially'' huge with $\omega$ supercompacts
above it, is somewhat technical, but it may be obtained with the following
more familiar assumptions:
\begin{quote}
A cardinal $\kappa$ and an increasing sequence $\la \lambda_i \mid i <
\omega \ra$ such that:
\begin{description}
\item[$B_1$:] For $\lambda=\sup \{ \lambda_i \mid i < \omega \}$ and $\mu =
\lambda^+$, each $\lambda_i$, $ i > 0$, is $\mu$-supercompact.
\item[$B_2$:]$\kappa$ is the critical point of an embeding $j : V
\longrightarrow M$ where $j(\kappa) = \lambda_0$ and $M^\mu \subseteq M$.
\end{description}
\end{quote}
Our aim in this subsection is to prove that if cardinals $\kappa$ and $\la
\lambda_i \mid i < \omega \ra$ satisfy $B_1$ and $B_2$, then there is a
generic extension in which cardinals that satisfy $A_1$ and $A_2$ can be
found.

Let $\rho$ be a cardinal, and $L=\seqlambda$ be
an increasing sequence
 with limit $\lambda$ of cardinals such
that, $\lambda_0=\rho$ and
for $i > 0$, $\lambda_i$ is $\lambda^+$-supercompact. Than $L$
 is called ``the
minimal supercompact
sequence above $\rho$'' if $\lambda$ is the least cardinal
above $\rho$ such that the interval $(\rho, \lambda)$
possesses an $\omega$-sequence of $\lambda^+$-supercompact cardinals,
$\lambda_0=\rho$,
 and  each
$\lambda_{i+1}$ is the first $\lambda^+$-supercompact cardinal above
$\lambda_i$.

Fix a function $g$ such that, for every $\rho$, $g(\rho)= \seqlambda$ is
such that $\lambda_0=\rho$ and $\la \lambda_i \mid  i < \omega \ra$
 is the
minimal supercompact sequence above $\lambda_0$ (if it exists, and $g$ is
undefined otherwise).

Looking at $j^2$, it is not difficult to see that if $\kappa$ is huge, then
it is $j(\kappa)$-supercompact, not only in $V$ but in $M$ as well (use the
combinatorial characterization of supercompactness). Hence, as high as we
wish below $\kappa$, there are cardinals that are $\kappa$-supercompact.
So, for every $\rho < \kappa$, $g(\rho) < \kappa$ (is defined), and
thus, for every $\rho < j(\kappa)$, $g(\rho) < j(\kappa)$ as well.

Let $\kappa$, $\mu$, and $\seqlambda$
be  as in $B_1$ and $B_2$. So $\kappa$ is
 the critical point of an embedding $j: V \longrightarrow M$, where $M^\mu
\subseteq M$, and (by taking the minimal sequence above $j(\kappa)$)
we may assume that $$g(j(\kappa))=\seqlambda=L.$$ Thus, if $L_0$ is the
minimal supercompact sequence above $\kappa$, and $\tau= (\sup(L_0))^+$,
then $j(L_0)=L$, and
 $j(\tau)=\mu$.
 $\kappa$ is thus $\tau$-huge (as in Lemma 2.3).

For every ordinal $\alpha \leq j(\kappa)$ a cardinal $\rho_0(\alpha)$
and an iteration $P_\alpha$ of length $\alpha$ with Easton support, is
defined by induction. Then, finally, $P_{j(\kappa)}$ is the required poset
which gives an extension where $\kappa$ and $L$ satisfy $A_1$ and
$A_2$.
\begin{enumerate}
\item
$P_0$ is the trivial poset and $\rho_0(0)=\aleph_1$.
\item For limit $\alpha$, $P_\alpha$ consists of all partial functions $f$
defined on $\alpha$ such that $f \rest \gamma \in P_\gamma$ for all $\gamma
< \alpha$, and $f$ has the Easton support property: $\dom(f) \cap \gamma$
is bounded below $\gamma$ for every inaccessible $\gamma$.
The cardinal $\rho_0(\alpha)$ is the
 first inaccessible cardinal above
all the $\rho_0(\gamma)$'s, $\gamma < \alpha$.
\item If $P_\alpha$ and $\rho_0(\alpha)$ are defined, then
$P_{\alpha+1}= P_\alpha \ast \Coll(L^\alpha)$ for $L^\alpha
=g(\rho_0(\alpha))$. The
first inaccessible above the cardinals in $L$ is $\rho_0(\alpha+1)$.
\end{enumerate}

Standard arguments prove that for Mahlo  $\gamma$'s that are closed
under the function $\alpha \mapsto \sup(g(\alpha))$, $P_\gamma$ satisfies
the $\gamma$-c.c. Also, for any $\alpha < \beta$, $P_\beta$ can be
decomposed as  $P_\alpha \ast R$, where $R$ is defined in $V^{P_\alpha}$
as an Easton support iteration of collapses determined by the same function
$g$, but beginning with $\rho_0(\alpha)$. (Supercompact cardinals remain
supercompact in any generic extension done via a poset of smaller size; see
LLL.)

Set $P=P_{j(\kappa)}$. We will show in $V^P$ that $\kappa$ and
$L=g(j(\kappa))$ satisfy the properties $A_1$ and $A_2$. Recall that
$\mu=\lambda^+$ where $\lambda=\sup(L)$. In $V$, we have an
elementary embedding $j:V \longrightarrow M$  into a transitive inner model
$M$ such that $M^\mu \subseteq M$. Again, the argument that small forcing
will not destroy supercompacness can show that the sumpercompact cardinals
in $L$ remain $\mu$-supercompact in $V^P$; that is, $A_1$ is easy.
  We promised to prove that in $(V^P)^{\Coll(L)}$,
$\kappa$ is $\tau$-huge for $\tau=\kpompo$, but in fact we will
 find a condition in $P\ast \Coll(L)$
and show that extensions through this
condition satisfy this requirement. The argument is fairly standard, but we
repeat it for completeness' sake.

Observe that
$$ j(P_{\kappa+1})= [ P_{j(\kappa+1)}]^M,$$
 but the
closure of $M$ under $\mu$ sequences implies that interpreting this
iteration in $V$ or in $M$ results in the same poset $Q$. So

 $$Q=  j(P_{\kappa+1})=
P_{j(\kappa)+1} = P_{j(\kappa)} \ast \Coll(L).$$
 We will find a
condition $q \in Q$ that forces $\kappa$ to be
$\tau$-huge (as in Lemma 2.3). In fact, the following suffices:
\begin{lemma}
\label{l24}
There is a condition  $q \in Q$
such that if $K$ is a  $V$-generic filter over $ Q$ containing $q$,
then the collection $\{ j(p) \mid p \in K \}$ has an upper bound in
$j(Q)/K$.
\end{lemma}
The meaning and proof of this lemma are clarified by the following:
Decompose
\begin{equation}
\label{A}
Q= P_{\kappa +1} \ast R
\end{equation}
where  $R$, the remainder, is an
Easton support iteration, starting above  $\Coll(g(\kappa))$,
of collapses guided by $g$,
going up to $j(\kappa)+1$.
Now apply $j$ to get
$$j( Q)=[P_{j(\kappa)+1} \ast j(R)]^M= [Q \ast j(R)]^M.$$
So $Q$ is a factor of $j(Q)$, and
$j(Q)/K$ can be formed in $M[K]$.
The lemma claims first that each
$j(p)$, for $p \in K$, is in $j(Q)/K$, and then that this collection has an
upper bound.

Since $2^\lambda = \mu$, the cardinality of $Q$ is $\mu$. It follows that
in $V[K]$ $(M[K])^\mu \subseteq M[K]$. Certainly, $j(Q)/K$ is
$\mu^+$-closed in $M[K]$ (in fact it is $\tau_0$-closed, where $\tau_0$ is
the first inaccessible in $M$ above $\lambda$). Hence $j(Q)/K$ is
$\mu^+$-closed in $V[K]$. So, to prove the lemma, we only need to choose $q
\in Q$ which forces that $$j(p) \in j(Q)/K$$ for $p \in K$.

Analyzing (\ref{A}), we write $p \in Q$ as $p=\la p_0,t,r \ra$ where $p_0
\in P_\kappa$, $t$ is a name, force to be in $\Coll(L)$, and $r$ is forced
to be in $R$. Then $j(p) = \la p_0, j(t), j(r) \ra$ (because $j(p_0)=p_0$
by the Easton condition). Now $j(p) \in j(Q)/K$ iff the projection of
$j(p)$ on $Q$, namely $\la p_0, j(t) \ra$, is in $K$.

The definition of $q$ can now be given, Define $q \in P_{j(\kappa)+1}$ as
$\la \emptyset, \sigma \ra$, where $\sigma \in V^{P_{j(\kappa)}}$ is forced
to be in $\Coll(L)$. It is easier to describe the interpretation of
$\sigma$ in $V[H]$, where $H \subset P_{j(\kappa)}$ is $V$-generic. Well,
look at all conditions $p \in P_{\kappa+1} \cap H$ (there are $< j(\kappa)$
of them); write each such $p$ as $p=\la p_0, t \ra$; interpret $t$ as a
condition in $\Coll (L)$, and take the supremum in $\Coll(L)$ of all of
these condition. This proves the lemma and we now see how the result
follows.

\begin{lemma}
Assuming $q$ is as in Lemma \ref{l24}, $q \forces_Q \kappa \mbox{\ is \ }
\tau-\mbox{huge}.2.3$.
\end{lemma}
{\bf Proof}. Let $K \subseteq Q$ be a $V$-generic filter containing $q$.
Work in $V[K]$ and let $s_0 \in j(Q)/K$ be an upper bound of $\{ j(p) \mid
p \in K \}$. We are going to define in $V[K]$ an ultrafilter ${\cal U }$
over $\powertaumu$ that satisfy the properties of Lemma 2.3. For this, we
fix $\la A_\xi \mid \xi \in \mu^+ \ra$, an enumeration af all subsets of
$\powertaumu$, and plan to decide inductively whether $A_\xi \in {\cal U}$
or not. Construct by induction an increasing sequence $\la s_\xi \mid \xi
< \mu^+ \ra$ of conditions in $j(Q)/K$ as follows:
\begin{enumerate}
\item At limit stages, $\delta < \mu^+$, use the $\mu^+$-completeness of
$j(Q)/K$ to find an upper bound to $\la s_\xi \mid \xi < \delta \ra$.
\item If $s_\xi$ is defined, pick for $A_\xi$ a name $a_\xi$ such that
$a_\xi[K]$, the interpretation of $a_\xi$ in $V[K]$, is $A_\xi$. Then
$j(a_\xi) \in M^{j(Q)}$, and we find an extension $s_{\xi+1}$ of $s_\xi$
that decides whether $(j''\tau)\in j(a_\xi)$ or not. If the decision is
positive, then put $A_\xi \in {\cal U}$, and otherwise not. Two comments
are in order for this definition to make sense:
\begin{enumerate}
\item First, $j(a_\xi)$ is not a name in $j(Q)/K$-forcing, but in $j(Q)$.
Yet, from $M[K]$ any generic extension via $j(Q)/K$ takes us into a universe
that is also a $j(Q)$ generic extension of $M$, and it is as such that we
ask about the interpretation of $j(a_\xi)$.
\item Apparently, this definition depends on a particular choice of a name
for $A_\xi$, but in fact if $a_\xi'$ is another name, then the same
answer is obtained. The point of the argument is that some condition $p$ in
$K$ forces $a_\xi = a_\xi'$, and hence $[j(p)\forces_{j(Q)} j(a_\xi)=
j(a_\xi')]^M$. But since $j(p) \in j(Q)/K$ is extended by $s_0$, it can be
seen that the answer to the $\xi$'s question does not depend on the
particular choice of the name.
\end{enumerate}
\end{enumerate}
We leave it to the reader to prove that ${\cal U}$ thus defined satisfies
the required properties of Lemma 2.3. For example, let us prove that ${\cal
U}$ is $\kappa$-small below $\tau =\kpompo$ in $V[K]$. For some $\xi \in
\mu^+$, $A_\xi = \{ A \in \powertaumu \mid | A \cap \tau | < \kappa \}$, and
a name $a_\xi$ for $A_\xi$ was chosen and an extension $s_{\xi+1}$ deciding
whether $j '' \mu \in j(a_\xi)$ was thought after. But some condition $p
\in K$ forces $a_\xi =
\{ A \in \powertaumu \mid | A \cap \tau | < \kappa \}$, and hence $j(p)$,
and decidedly $s_0$, forces $j(a_\xi)= \{ A \in P^{j(\tau)}_{j(\mu)} \mid | A
\cap j(\tau) | < j(\kappa) \}$. Now $j'' \mu = A$ has order-type
$\mu=j(\tau)$, and $A \cap j(\tau) = j'' \tau$ has cardinality $\tau <
j(\kappa)$.
\section{There are no  Aronszajn trees on successors of singular limits of
compact cardinals}
\label{s3}
The paper really begins here, with the following theorem of ZFC.
\begin{theorem}
If $\gl$ is  singular and a limit of strongly compact
cardinals, then there are no $\gl^+$-\A\ trees.
\end{theorem}
{\bf Proof}.  For notational simplicity,  assume that $\cf(\gl)=\go$.
Let $\la \gl_i \mid i < \go \ra$ be an increasing $\go$-sequence of
strongly compact cardinals with limit $\gl$.  (Recall that a cardinal
$\gk$ is strongly compact if every $\gk$-complete filter can be
extended to a $\gk$-complete ultrafilter.)
Let $T$ be a tree of height $\gl^+$ and levels of size $\leq \lambda$
(a $\gl^+$--tree), and we will find a $\gl^+$ branch in $T$.  We
may assume that $T_{\ga}$, the $\ga$th level of $T$, is the set $\gl
\times \{\ga \}$.  Accordingly, we define $T_{\ga,n}= \gl_n \times \{\ga
\}$, so that $T_{\ga} = \bigcup_{n < \go} T_{\ga,n}$.
 The proof for the existence of the  branch is divided into two
steps:
\begin{description}
\item[Step one:] There is an unbounded $D \subseteq \gl^+$ and a fixed
$n \in \go$ such that
whenever $ \ga < \gb$ are both in $ D$, then, for some $a \in
T_{\ga, n}$ and $ b \in T_{\gb ,n },\ a <_T b$.
We call an unbounded set $D$ and a collection $\la T_{\ga, n } \mid
\ga \in D \ra$ as above a {\em spine} of $T$.  Thus the the first part of the
proof provides a spine for every $\gl^+$ tree.

Indeed, using the fact that $\gl_0$ is strongly compact, extend the
co--bounded subsets of $T$ (that is, those
subsets whose complement has
 cardinality $\leq \lambda$) to a countably complete, uniform
 ultrafilter, $u$, over
$T$.  Given $\ga \in \gl^+$ (considered as a level of $T$)
   define $n_{\ga} \in \go$ by the
following procedure:  For every $x \in T$ of level $> \ga$, let $r^\ga_x \in
T_{\ga}$ be such that $r^\ga_x <_T x$, and set $n=n_x$ to be (the least) such
that $r^\ga_x \in T_{\ga,n}$.  Since the set $T \setminus (T \rest \ga
+1)$ is
in $u$, it follows from the $\aleph_1$--completeness of $u$ that for some
$n=n_{\ga},\ \{ x \in T \mid n_x \leq n \}=X_{\ga} \in u$.

Now there is an unbounded $D \subseteq \gl^+$ and a fixed $n$ such
that $n=n_{\ga}$ for $\ga \in D$. If we take any two ordinals
 $\ga <  \gb$ in $D$,
 then the intersction $X_{\ga} \cap X_{\gb}$ is in $u$, and any $x$ in this
intersection is such that $a=r^\ga_x$ and $b=r^\gb_x$ are comparable
(being both below $x$),
and in the $n$th part, as required.
\item[Step two:] Every spine has a cofinal branch.  Suppose that
 $D$ and $n$ define
 a spine of $T$ as above. That is, assume $\{ T_{\ga,n} \mid \ga \in D \}$,
where $D \subseteq \lambda^+$ is unbounded, is a collection such that for
every $\ga < \gb$ in $D$ there are $a \in T_{\ga,n}$, $b \in T_{\gb,n}$ such
that $a <_T b$.  Find a $\gl_{n+1}$--complete ultrafilter,
 $v$, over $\gl^+$
containing $D$ and the co-bounded subsets.  Fix any $\ga \in D$. For every
$\gb > \ga$ in $D$ find $a(\gb) \in T_{\ga, n}$ and $b(\gb) \in
T_{\gb,n}$ such that $a(\gb) <_T b(\gb)$. Use the completeness of $v$,
and the fact that the cardinality of each level of the spine is only
$\lambda_n$, to find
$a_{\ga} \in T_{\ga,n}$ and $\xi_{\ga} \in \gl_n$ such that	for a set of
$\gb$'s in $v$, $a_{\ga}= a(\gb)$ and $b(\gb)= \la \xi_\ga , \gb \ra$
(which is
the $\xi_\ga$--th element of $T_{\gb,n}= \lambda_n \times \{ \gb \}$).
 For an unbounded $D' \subseteq
D$ the ordinal $\xi_\ga$ has the fixed value $\xi$ for $\ga \in D'$. Now
the collection $\{ a_{\ga} \mid \ga \in D' \}$ is a branch of $T$.
\end{description}

In a very direct way, one can generalize this to find that if $\gl$ is
singular and a limit of strongly compact cardinals as above, then any
 strong system over $\gl^+$ with index-set of size $< \gl$ has a
branch of size $\gl^+$.

It may be illuminating to remark that the structure of this theorem
resembles somewhat the structure of the consistency proof for $\alephomp$.
Assuming that $\lambda$ is a singular limit of strongly compacts, we have
proven in Step one that any $\lambda^+$-tree (viewed as a strong system)
posseses a narrow subsystem. Then, in Step two, it was shown in fact that
every narrow system (with fewer than $\lambda$ relations) has a
 $\lambda^+$ branch. Now, when proving the consistency theorem (assuming
some large cardinals), we will find in the first step a generic extension
such that the pre-system of every strong system on $\alephomp$, with
countably many relations, has a narrow subsystem (that is, one of width
$\aleph_n$) in the ground model (Theorem  ).
 So, if in the ground model every narrow
system has a cofinal
branch, then we are done. Unfortunately, we only know how to
arrange the ``potential branching property'' in the ground model: Every
narrow system has a branch in some sufficiently closed generic extension.
It turns out that this is sufficient, because of a ``preservation theorem''
(Theorem \ref{preservation}).

\section{The narrowing property}
\label{s4}
\begin{definition}
Let $ \mu > \chi$ be two cardinals, and $\mu=\lambda^+$.
A poset $Q$ has the narrow pre-system property for
 $(\mu,\chi)$ if whenever $\cal S$ is a $Q$-name of a strong $\mu$-system,
with $\leq \chi$ relations, then $\Pre_Q(\CS)$ has a narrow subsystem.
\end{definition}
\begin{theorem}
\label{narrowing}
  Suppose that $\gk$ is $\kpompo$-huge. That is, $\kappa$
 is the critical point of an elementary embedding $j
: V \rightarrow M$, where $M$ is a transitive class such that $M^\mu
\subseteq M$, for $\mu = \gl^+= \jkpompo$.
Let $Q=\Coll(\kpompo, < j(\gk))$.  Then  $Q$ has
the  narrow pre-system property for $(\mu, \kpom)$.
\end{theorem}
{\bf Proof}.
Suppose $\CS_0$ is forced by every condition in $Q$ to be a strong system
on $\mu$, and $\kappa^{+\omega}$ is its index set. Let $\CS_1=\Pre_Q(\CS_0)$ be
its pre-system; we must find a narrow subsystem of $\CS_1$ (that is, one of
width and index set $< \lambda$).
The $\alpha$th level of $\CS_0$ (and $\CS_1$) is $\lambda \times \{ \alpha
\}$, and we denote it by $(S_0)_\alpha$. The $n$th part of this level,
$j(\kappa)^{+n}$ is denoted $(S_0)_{\alpha,n}$.

Observe that $j(Q)$ is $\mu$--closed in $M$ since $Q$ is
$\kpompo$--closed and $\mu=j(\kpompo)$. In fact, $j(Q)$ is $\mu$--closed
in $V$ since $M$ is sufficiently closed.

$j(\CS_0)$ is in $V^{j(Q)}$ a strong $j(\mu)$-system with relations indexed
by $j(\kappa)^{+\omega}=\lambda$. The $\alpha$th level of $j(\CS_0)$ is
$j(\lambda) \times \{ \alpha \}$. It is more convenient to denote this level
by $(jS_0)_\alpha$. Similarly, the $n$th part of this level is denoted
$(jS_0)_{\alpha,n}$ ($= j(\kappa)^{+n} \times \{\alpha\}$).

It follows from the closure of $M$ under $\mu$-sequences
 that $j''\mu$ is a bounded
subset of $j(\mu)$ in $M$,
 and we let $\gb^\ast < j(\mu)$ be a bound of $j''\mu$.
Let $b^\ast$ be any fixed ordinal in $j(\lambda) \times \{ \beta^\ast \}$
(so $b^\ast$ is any node of level $\beta^\ast$ in $j(\CS_0)$).

Inductively,
define---in $M$---a $j(Q)$ increasing
sequence of conditions $\{ s_{\ga} \mid \ga < \mu\} $, starting with any
condition, as follows:
\begin{enumerate}
\item At limit stages, the $\mu$-closure of $j(Q)$ is used to find an upper
bound to the sequence of length $< \mu$ so far constructed.
\item If $s_\alpha$ is defined, then $s_{\alpha+1}$ is defined as follows:
Since $j(\CS_0)$ is forced to be a strong system, there
exist $a \in (jS_0)_{j(\alpha)}$ and $\zeta < \lambda$ such that $\la a , b^\ast \ra$ is forced by some
extension of $s_\alpha$ to stand in the $\zeta$th relation of $j(\CS_0)$.
So we pick $s_{\alpha+1}$, extending $s_\alpha$, $a_\alpha \in (jS_0)_{j(\alpha)}$,
and $\zeta_\alpha$ such that $$s_{\alpha+1} \forces_{j(Q)} \la a_\alpha, b^\ast
\ra\ \mbox{stands in the}\ \zeta_\alpha \mbox{th relation}.$$
\end{enumerate}
Since $\lambda < \mu$ there is a fixed $\gz^0 < \lambda$, and $n \in \omega$,
such that for some  unbounded
set $D \subseteq \mu$, $\gz_\ga = \gz^0$ and $a_\alpha \in (jS_0)_{\alpha, n }$
 for all $\ga \in D$.

Now the pre-system $\CS_1$ has width $\lambda$ and relations indexed by
$ Q\times \kpom$. We claim that the narrow substructure of $\CS_1$
defined by $\{ (S_0)_{\alpha,n} \mid \alpha \in D \}$ is a system, thereby
proving the theorem. If this is not the case, then for some $\alpha_1 <
\alpha_2$ in $D$, ``there are no $a_1 \in (S_0)_{\alpha_1,n}$ and
$a_2 \in (S_0)_{\alpha_2,n}$ such that $\la a_1,a_2 \ra$ stands in a
relation of $\CS_1$'', or specifically, ``there are no  $a_1 \in
(S_0)_{\alpha_1,n}$, $a_2 \in (S_0)_{\alpha_2,n}$, $q \in Q$, and $\zeta
\in \kpom$ such that $q \forces_Q \la a_1, a_2 \ra $ stand in the
$\zeta$th relation.'' But then, applying $j$ to this statement
we get a contradiction to: $s_{\alpha_2 + 1} \forces
\la a_{\alpha_1}, b^\ast \ra,\
\la a_{\alpha_2}, b^\ast \ra\ \mbox{and hence}\ \la a_{\alpha_1},
a_{\alpha_1} \ra \ \mbox{as well}\ \mbox{stand in the}\ \zeta^0
\mbox{th relation}.$

\begin{cor}
\label{c42}
Let $\kappa$ be $\kpompo$-huge (as in Lemma 2.3). Then  it is possible to
collapse $\kpompo$ to be $\alephomp$ with a forcing poset that has the
narrow pre-system property for $(\jkpompo, \omega)$. In other words,
there is a forcing
poset $P$ such that
\begin{enumerate}
\item $j(\kappa)^{+(\omega+1)}$ becomes $\alephomp$ in $V^P$, and
\item the pre-system$_P$ of every strong system on $\jkpompo$
with countably many
relations in $V^P$ has a narrow subsystem.
\end{enumerate}
\end{cor}
{\bf Proof}.
The desired poset is simply the collapse of $j(\kappa)$ to become
$\aleph_2$, but not in the most direct way.
It is rather the product of two collapses that works: the collapse of
$\kpom$ to $\aleph_0$, and the one that makes $j(\kappa)$ the double
successor of $\kpom$ (both posets are defined in $V$).
Let $Q=Coll(\kpompo, < j(\gk))$,  $K=  Coll(\aleph_0,\gk^{+ \omega})$,
and then $P=Q \times K$ is the desired collapse.
In $V^P$, $\kpom$ is countable, $\kpompo$ is $\aleph_1$,
$j(\gk)$ is $\aleph_2$, and $\mu=\lambda^+=\jkpompo$ becomes $\aleph_\omega^+$.

So let $\CS_0$ be in $V^P$ any strong $\mu$-system; then $\Pre_P(\CS_0)$
can be obtained in two stages, corresponding to the product $P=Q \times K$
and decomposition $V^P=(V^Q)^K$. First, in $V^Q$, form $\CS_1=\Pre_K(\CS_0)$. Then
$\CS_1$ is in $V^Q$ a strong system on $\mu$, with $|K| \times \aleph_0=
\kpom$ relations, one relation for each pair formed with a condition in
$K$ and a relation (index) in $\CS_0$. Hence, by the theorem,
$\Pre_Q(\CS_1)\sim \Pre_P(\CS_0)$ has (in $V$) a narrow subsystem as
required.

\section{The potential branching property and a model with no Aronszajn
trees}
\label{s5}
Suppose $\mu=\lambda^+$.  The potential branching property for $\mu$ is
the following statement:
\begin{quote}
If $\CS$ is a narrow system on $\mu$, then for every $\chi < \lambda$
there is a $\chi$-complete forcing poset that introduces an unbounded
branch to $\CS$.
\end{quote}
Recall that a branch of a system is a set of nodes and a relation in the
system which includes every pair from the set.
In the following section we will see how to obtain the potential
branching property, but here we use it to obtain a model with no
Aronszajn trees on $\alephomp$.
\begin{theorem}
\label{notrees}
  Let $\gk$ be $\kpompo$-huge (as in Lemma 2.3).
  Suppose that the potential branching property for
$\mu$ holds ($\mu= j(\kappa)^{+(\omega+1)}$).
Then there is a generic extension in which $\mu$ becomes
$\aleph_\omega^+$ and it carries no
\A\ trees.
\end{theorem}
{\bf Proof.}
The poset $P$ of Corollary \ref{c42}  works. Recall that
$Q=Coll(\kpompo, < j(\gk))$,  $K=  Coll(\aleph_0,\gk^{+ \omega})$,
and then $P=Q \times K$ is the desired collapse.
We will prove
that there are no $\mu$--\A\ trees in $V^P$.
  The demonstration of this result
depends on a preservation theorem which will only be proved in the following
subsection, and it goes as follows.

 Suppose that
 $\tildaT$ is in
$V^P$ a name of a  $\mu$-tree.  We will first show that there is in
$V$ a $\kpompo$--complete poset $R$ such that, in $(V^P)^R$, $\tildaT$
aquires an unbounded branch. Then
the preservation theorem (\ref{preservation}, applied with
$\lambda=j(\kappa)^{+\omega}$)
 shows that  $\tildaT$
has a branch already in $V^P$.

 To see how $R$ is obtained,  apply the corollary to $\tildaT$, considered
as a one relation strong system, and find a narrow subsystem $\CS$ to
$\Pre_P(\tildaT)$. But then for $\chi=|P|<|lambda$ there is a
$\chi^+$-complete forcing poset $R$ that introduces an unbounded branch to
$\CS$. This can be shown to give an unbounded branch to the tree $\tildaT$
in $V^{P\times R}$. But then the following preservation theorem shows that
$\tildaT$ already has a branch in $V^P$.

\subsection{The preservation theorem}
\begin{theorem}
\label{preservation}
Let $\lambda$ be a singular cardinal (of cofinality $\omega$, for
notational simplicity), and suppose that $P$ and $R$ are two posets such
that:
\begin{enumerate}
\item
$\|P\| = \chi < \lambda$,
and  $\tildaT$ is a $\lambda^+$-tree in $V^P$.
\item $R$ is $\chi^+$--closed.
\end{enumerate}
Then any $\lambda^+$-branch of $\tildaT$ in
 $V^{P \times R}$ is already in
 $V^P$.
\end{theorem}
{\bf Proof}.  As before, $T_\alpha = \lambda \times \{\alpha \}$ is the
$\alpha$th level of $\tildaT$, for $\alpha < \mu= \lambda^+$.
 Let $\bfB \in V^{P \times R}$ be a name of a cofinal branch of
$\tildaT$, supposedly not in $V^P$. We also view $\bfB$
as a name in $(V^P)^R$ (that is, a name in $R$--forcing, in $V^P$).

Say (in $V^P$) that two conditions $r_1, r_2 \in R$ force
distinct values for $\bfB \cap T_\ga$ iff for some $a_1 \not =
a_2$ in $T_\alpha$,  $r_i \Rforces a_i \in \bfB$, for  $i=1,2$.
A weaker property, which even may hold when $r_1$ or $r_2$ do not determine
$\bfB \cap T_\ga$, is that whenever $r'_1$ and $r'_2$ are extensions of $r_1$
and $r_2$ that determine the value of $\bfB \cap T_\ga$ then $r'_1$ and
$r'_2$ force distinct values for $\bfB \cap T_\ga$. In this case we say
that $r_1$ and $r_2$ ``force contradictory information on $\bfB \cap
T_\ga$''.
 Observe that if $ \ga < \gb < \mu$ and $r_1,\ r_2$ force contradictory information
on $\bfB \cap T_\ga$, then they force contradictory information on
 $\bfB \cap T_\gb$.

Working in $V$,
our aim is to tag the nodes of the tree $\ltree$ (finite sequences from
$\lambda$) with conditions in $R$ that are forced to force
 pairwise contradictory
information on $\bfB \cap T_\ga$ for some $\ga$, and this will be shown to be
 impossible
because the cardinality of $T_\ga$ is $\lambda$ and $\ltree$ has $\geq \gl^+$
branches.  We will denote the tag
of $\gs \in \ltree$ with $r_\gs \in R$.  The required properties of this
tagging are the following.
\begin{enumerate}
\item If $\gs_1 \subset \gs_2$ in $\ltree$, then $r_{\gs_1} \leq
r_{\gs_2}$ in $R$.
\item For every node $\gs \in \ltree$, for any two immediate extensions
$\gs_1, \gs_2$ of $\gs$, there are an ordinal $\ga= \ga(\gs_1, \gs_2)$,
 and a dense set $D
\subseteq P$, such that for every $p \in D$
\[ p \Pforces \ r_{\gs_1} \mbox{\it  and } r_{\gs_2}\ \mbox{\it  force
contradictory information on } \bfB \cap T_\ga.\]
\end{enumerate}
Why this suffices?  Because, assuming such a construction, let $\gb <
\mu$ be above all the ordinals $\ga(\gs_1,\gs_2)$
 mentioned in item 2 for nodes
$\gs_1$ and $\gs_2$, and look at the set of all full branches
$\lambda^\omega$. For each $f \in \lambda^\omega$, let $r_f\in R$  be
an extension of the
upper bound of the conditions $r_{f \rest n}$, tagged along the branch
$f$, an extension that determines the value of $\bfB \cap T_\gb$.
  We claim now that if $f \not = g$ are full branches, then there is
a dense set $D \subseteq P$ such that for every $p \in D$, $p \Pforces\
r_f \ \mbox{\it and }\ r_g\ \mbox{\it force distinct values for}
\ \bfB \cap T_\gb.$ Indeed, let $\sigma \subset f \cap g$ be the
splitting node, then item 2 gives the required dense set.  To conclude
the proof, we find that, in $V^P$, any two branches of $\lambda^\omega$
give distinct values for $T_\gb$, and since $\lambda^{\aleph_0}
\geq \lambda^+$, this shows that $T_\gb=\lambda \times \{ \gb \}$
 contains $\lambda^+$ distinct
nodes in $V^P$, which is not possible since $\gl^+$ is not collapsed.

In the construction of the tags, the following lemma, stated in general
for any poset $S$,
provides a basic tool.
\begin{lemma}
\label{wide}
Let $S$ be any forcing poset. Suppose that $\lambda$ and $\lambda^+$ are
 a cardinal and its successor,  $T$ is a $\lambda^+$--tree,
and $\bfB$ is a name of a  $\lambda^+$--branch of $T$ in $V^S$.
   If $\bfB$ is a new branch ($\bfB$ is not in
$V$), then for some $\ga$ there is a set $X \subset	T_\ga$ of
cardinality $\lambda$, in $V$, such
that every $x \in X$ is forced by some condition in $S$ to be in $\bfB$.
\end{lemma}
{\bf Proof}.
We will say that $s \in S$ is $\lambda$--wide at $T_\ga$ if there are
$\lambda$ extensions of $s$ that force pairwise distinct values for $\bfB
\cap T_\ga$.  If we start with an arbitrary condition, our proof will
give that every $s \in S$ is $\lambda$--wide
at some $T_\ga$, $\ga < \lambda^+$.
	
 Define $E$, in $V$, to be
 the set of possible nodes of $\bfB$: $$E = \{ a
\mid \ \mbox{some condition in}\ S\ \mbox{forces}\ a \in \bfB \}.$$  We
want some $\ga < \lambda^+$ such that $\mid E \cap T_\ga \mid =
\lambda$.  So assume, on the contrary,
 that $\mid E \cap T_\ga \mid < \lambda$ for every
$\ga$. Then  $E \subseteq T$ satisfies the following properties:
\begin{enumerate}
\item Any node in $E$ has extensions in $E$ at arbitrarily higher level.
\item $E$ is downward closed in $T$.
\item Any node in $E$ has two incomparable extensions in $E$ (for
otherwise, a condition would force that $\bfB$ is
 in $V$).
\item For every $\ga< \lambda^+$, $\mid E \cap T_\ga \mid < \lambda$.
\end{enumerate}
This is not possible: let $U \subseteq \lambda^+$ be a closed unbounded
set such that if $\ga \in U$ then whenever $\gamma < \ga$, if
$a \in E\cap T_\gamma$,
then $a$ has two incomparable extensions of height $< \ga$ in $E$.
 Then pick
any $\ga \in U$ such that $\ga \cap U$ has order--type $\geq \lambda$,
and conclude that $E \cap T_\alpha$ has size $\lambda$ by splitting the points
at levels in $U$ below $a$. QED

Returning to our case $S=P \times R$, we will say that $(p,r) \in P
\times R$ is $\lambda$--wide at $T_\ga$ if there are $\lambda$
extensions of $(p,r)$ that force pairwise distinct values for
$\bfB \cap T_{\ga}$. (This is not exactly the same definition as the one
given above, because $\tildaT$ is not assumed to be in $V$; however, its
level--sets are, and so this definition is meaningful.)

\begin{lemma}
\begin{enumerate}
\item Any condition $(p,r) \in P \times
R$ is $\lambda$--wide at some $T_\ga$.
\item If $(p,r)$ is $\lambda$-wide at $T_\ga$,
then it is also $\lambda$--wide at any higher level $T_\gb$.
\end{enumerate}
\end{lemma}
{\bf Proof}.
 Indeed, let $G$ be a $V$--generic
filter over $P$ containing $p$.  In $V[G]$, $\lambda$ and $\lambda^+$
are not collapsed, and we can use Lemma \ref{wide}
  to find some $\ga < \mu$ such
that there are $\lambda$ possible values for $\bfB \cap T_\ga$ (forced by
some extensions of $r$). Any such value is also a possible value for some
extension of $(p,r)$ (in $V$), and hence when the set of possible values
for $\bfB \cap T_\ga$ is calculated in $V$ it must have cardinality
$\lambda$ as well.

  For the second part assume that $(p,r)$ is $\lambda$ wide at $\ga$.
  given any $\gb > \ga$, find first  $(p_i, r_i)$ extending $(p,r)$, for $i <
\lambda$, that determine distinct values of $\bfB \cap T_\ga$
and then extend each pair to a condition $(p'_i,r'_i)$ that determines $\bfB
\cap T_\gb$.  Even though it may be possible for two such extensions to
determine the same point in $T_\gb = \lambda \times \{ \gb \}$,
 it is not possible for $\| P \|^+$
extensions to determine the same point (because in such a case we would
have two extensions with compatible $P$ coordinates, and this is not possible).
Hence, as $\| P \| < \lambda$,
 there are $\lambda$ possible values for $\bfB\cap T_\gb$.
Observe, however, that if $(p,r)$ is $\lambda$--wide at $T_\ga$, then
extensions of $(p,r)$ need not be $\lambda$--wide at the same $T_\ga$, and it
may be necessary to go to higher levels.

\begin{lemma}
If $\{ r_j \mid j < \lambda \} \subseteq R$, and $p_0 \in P$ are given,
then, for some ordinal $\ga$,
 there are extensions $r'_j \geq r_j$ in $R$, for $j < \lambda$,
 such that for every pair $i
< j$ there is $p_1 \geq p_0$ in $P$  such that $p_1
\Pforces\ r'_i \ \mbox{and}\ r'_j\ \mbox{ force distinct values for }
 \bfB\cap T_\ga$.
\end{lemma}
{\bf Proof}.
First, by our last lemma,
  find for every $i<\lambda$ an ordinal
$\ga_i$ such that $(p_0,r_i)$ is $\lambda$--wide at $T_{\ga_i}$,
and then
let $\ga$ be above all of these $\ga_i$'s.
By the second part of the lemma,
 each $(p_0,r_i)$ is $\lambda$--wide at $T_\ga$. Now, by induction on $i <
\lambda$, we will define an extension $r_i' \geq r_i$, and two functions,
$e_i$ and $f_i$,  $e_i:P \rightarrow P$, and $f_i:P \rightarrow T_\ga$, such
that:
\begin{enumerate}
\item For every $a \in P$,  $e_i(a)$ extends $a$, and
$ (e_i(a), r_i') \Vdash_{P\times R}\ \bfB \cap
T_\ga = \{ f_i(a) \}$.
\item If $k<i<\lambda$ then
$$ f_i(p_0) \not \in \{ f_k(a) \mid a \in P \}.$$
That is, the value of $\bfB \cap T_\ga$ that $(e_i(p_0),r_i')$ determines
is different from all the values determined by previous conditions.
\end{enumerate}
Suppose that it is the turn of $r_i'$, $e_i$, $f_i$ to be defined.
Let $P= \{ p(\xi) \mid \xi < \chi \}$ be an enumeration of $P$, starting
with the given condition $p(0)=p_0$.
 By	induction on $\xi < \chi$, we shall define a condition $r_i^\xi \in
R$, and the values
 $e_i(p(\xi)) > p(\xi)$, and $f_i(p(\xi)) \in T_\ga$ such that:
\begin{enumerate}
 \item $\la r_i^\xi \in R \mid \xi < \chi \ra$ form an increasing sequence
 of conditions extending $r_i$
\item 	$(e_i(p(\xi)),r_i^\xi) \forces \bfB \cap T_\ga= \{
f_i(p(\xi))\}$.

\item $f_i(p(0)) \not \in \{ f_k(a) \mid a \in P,\; k<i \}$.
\end{enumerate}
First,  use the fact that $(p_0,r_i)$
is $\lambda$--wide at $T_\ga$
to find an extension $(e_i(p_0),r_i^0) \geq (p_0,r_i)$
that forces $\bfB \cap T_\ga = \{ f_i(p_0) \}$ for a value $f_i(p_o)$ that
satisfies (3) above.  Then construct the increasing sequence $r_i^\xi$ and
the values of $e_i$ and $f_i$ (using the $\chi^+$ completeness of $R$ at
limit stages), and finally define $r_i'$ to be an
upper bound in $R$ of that sequence.

Let us check that the requirements of the lemma are satisfied for $r'_i$.
If $k<i$
is any index, look at $p'=e_i(p_0)$, and let $p_1=
e_k(p')$.  Then $p_1$
is as required, because the value of $\bfB \cap T_\ga$ determined by
 $(p_1,r_i')$ (namely $f_i(p_0)$)
is distinct from the one determined by $(p_1,r_k')$ (namely $f_k(p')$).
This proves the lemma, and the following completes the proof of the
theorem  by showing how the tagging can be done.
\begin{lemma}
If $r \in R$, then there are extensions
$r'_i \geq r$ for $i< \lambda$ such that, for some
$\ga$, if $i < j< \lambda$ then  for some dense set $D=D_{i,j} \subseteq P$,
for every $p \in D$,
\[ p \Pforces r'_i \mbox{ and } r'_j \mbox{ force contradictory
information on } B \cap T_\ga. \]
\end{lemma}
{\bf Proof}. Enumerate $P = \{ p(\xi) \mid \xi < \chi \}$.
  By induction on $\xi \leq \chi$ we define:
\begin{enumerate}
\item  A sequence of conditions in $R$, $\la r_i^\xi \mid i < \lambda \ra$.
\item A family $D^\xi(i,j) \subset P$,
 increasing with $\xi$, for every $i < j < \lambda$. Finally, we will set
$D(i,j) = D^\chi(i,j)$, and to ensure that $D(i,j)$ is dense we demand that
$p(\xi)$ has an extension in $D^{\xi+1}(i,j)$.
\item An ordinal $\ga(\xi) < \mu$.
\end{enumerate}
We require that
for each $i$, $\la r_i ^\xi \mid \xi < \chi \ra$ forms an
increasing sequence, beginning with $r^0_i=r$.
 (Finally, $r'_i = r^\chi _i$ will be the required extension.)

 At limit stages $\gd$, $r^\gd_i$ is
an upper bound in $R$ of the conditions $r_i^\xi$, $\xi < \delta$.
$D^\delta(i,j)$ is the union of $D^\xi(i,j)$ for $\xi < \delta$.

At successor stages, $\xi+1$, the extensions
$\{ r_i^{\xi+1} \mid i< \lambda \}$ are defined using
the previous lemma for the collection $\{r_i^\xi \mid i < \lambda \}
$ and the condition $p_0=p(\xi)$. The lemma gives
 an ordinal $\ga=\ga(\xi)$ and extensions $p_1(i,j) \geq p(\xi)$,
such that $$ p_1(i,j) \forces r_i^{\xi+1}\ \mbox{\it and}\ r_j^{\xi+1}\
\mbox{\it force distinct values for}\ \bfB \cap T_\ga.$$
Then we define $D^{\xi+1}(i,j)$ by
$D^{\xi+1}(i,j)= D^\xi(i,j) \cup \{ p_1(i,j)\}$.
Finally, define $r'_i=r_i^\chi$, $\ga= \sup \{ \ga(\xi) \mid \xi < \chi \}$, and
$D_{i,j}= D^\chi(i,j)$. $D_{i,j}$ is dense in $P$, because every $p(\xi)$
has some extension in $D^{\xi+1}(i,j)$.

\section{The final model}
\label{s6}
For the consistency of {\em no Aronszajn trees on} $\alephomp$ we must show
how to obtain the assumptions of Theorem \ref{notrees}, namely how to get a
cardinal $\kappa$ which is $\kpompo$-huge
with the potential branching property for
$\jkpompo$.  We are going to do this in two stages: In the first, it is
shown that whenever an $\omega$ sequence of supercompact cardinals
conveging to $\lambda$ is collapsed, then the potential branching
property for $\lambda^+$ holds.  In the second stage, this is combined
with preparations described in Section \ref{s22}.

\begin{theorem}
\label{t61}
Suppose	 $\seqlambda$, with $\gl= \bigcup_{i<\go} \gl_i$ and $\mu=\gl^+$,
is an increasing $\go$--sequence of $\mu$--supercompact cardinals
(except for $\gl_0$ which is just a regular cardinal).
  Let $C=\Coll(\seqlambda)$
be the full support iteration that makes $\lambda_i$ to be
$\lambda_0^{+i}$. Then, in $V^C$, the potential
branching property holds for $\mu$:
\begin{quote}
 If $\CS$ is a narrow system on $\mu$,
  then, for every $k < \omega$, there is a $\gl_k$--complete forcing that
introduces an unbounded branch to $\CS$.
\end{quote}
\end{theorem}
{\bf Proof}.  We will actually prove the following combinatorial
statement in $V^C$:
\begin{quotation}
\noindent
 For every $n < \go$ and function $F: [\mu]^2
\rightarrow \chi$, where $\chi < \gl_n$, there is a
$\gl_n$--complete forcing $C^*$ such that in $(V^C)^{C^*}$ the following
holds: For some $\nu < \chi$ there is an unbounded set $U \subseteq \mu$
such that for every $\ga_1 < \ga_2$ in $U$ there is $\gb > \ga_2$ such
that \begin{equation}
\label{st}
F(\ga_1,\gb)=F(\ga_2,\gb)=\nu.
\end{equation}
\end{quotation}
  (We call such a set $U$ ``a branch'' of $F$.)

First, we shall
 see why this combinatorial statement suffices to prove the theorem.
Let $\CS=(T,R) \in V^C$ be a narrow system over $\mu$, and let
$\chi < \gl$ be such that the width of $\CS$ and the cardinality of its
index set are $\leq \chi$.  Suppose that $\chi < \gl_n$ and we will find
a $\gl_n$--complete forcing that introduces an unbounded branch to $\CS$.
For this,  define in $V^C$
a function $F: [\mu]^2 \rightarrow \chi^3$ by

\begin{quote}
 $F(\ga_1,\ga_2)=(\gz,\gt_1,\gt_2)$
iff the $\gt_1$ member of $T_{\ga_1}$ and the $\gt_2$ member of
$T_{\ga_2}$ stand in the $\gz$ relation $R_{\gz}$.
\end{quote}
 Then, by the assumed
combinatorial principle, there is an unbounded set $U \subseteq \mu$,
and fixed ordinals $\nu=(\gz, \gt_1,\gt_2)$ as in equation (\ref{st}).  This
 implies that the
$\gt_1$-th points of $T_{\ga}$ (that is $\la \gt_1, \ga \ra$) for $\ga \in U$
form an $R_{\gz}$ branch of $\CS$ (use item {\it 2} in the definition of
systems).

We now argue why it suffices to prove the combinatorial principle for
$n=0$.  Given $F: [\mu]^2 \rightarrow \chi$ where $\chi < \lambda_n$ in
$V^C$ (suppose for simplicity that {\em every} condition in $C$ forces
that $F$ is into $\chi$),  decompose $C \simeq C_n \ast
C^n$ where $C_n = \Coll(\la \gl_i \mid i \leq n \ra )$, and $C^n$ is the
name in $V^{C_n}$ of $\Coll(\la \gl_i \mid n \leq i < \go \ra )$.
In $V^{C_n}$ define $\gl'_m = \gl_{n+m}$.  Then each $\gl'_m$, for
$m>0$ is $\mu$--supercompact.  (Indeed the embedding $j:V\rightarrow M$
with critical point $\gl_k$, for $k > n$, can be extended in
$V^{C_n}$ to an embedding of $V^{C_n}$ into $M^{C_n}$,
 where $M^{C_n}$ posseses the
same $\mu$--closure properties.)
Thus, if we know case $n=0$ of the theorem in $V_0=V[C_n]$, we could apply it there to
$C^n=\Coll(\la \gl_i' \mid i < \omega \ra )$ and get in $V_0^{C^n} =V^C$ the
desired $\gl_0'$--complete ($\gl_n$--complete) poset that adds a branch
to $F$. To save ourselves from too many
superscripts, we denote $V^{C_n}$ by $V$ and $M^{C_n}$ by $M$ and
assume $n=0$.

So, returning to the theorem, assume that $G$ is a $V$-generic filter over
$C$, and $F$ is in $V[G]$ a function from $[\mu]^2$
into $\chi < \gl_0$. In the following lemma,
we will describe a $\lambda_0$-complete forcing $P$ in
$V[G]$ that introduces a $\mu$-branch to $F$.
 Let $j:V\rightarrow M$ be an
elementary embedding with critical point $\gl_1$, such that $j(\gl_1) >
\mu$ and  $M$ is
closed in $V$ under $\mu$-sequences.
The following lemma will be proved later on.
\begin{lemma}
\label{extension}
There is in $V[G]$ a $\lambda_0$-complete poset $P$ such that in $V[G]^P$
there is an extension of the embedding $j$ to an elementary embedding of
$V[G]$ into $N=M[j(G)]$.
\end{lemma}
 Accepting the lemma for a moment
$j(F) \in N$ can be defined; a function on $[j(\mu)]^2$ and into
 $\chi < \gl_0$ ($j$ is the identity below $\gl_1$).

 Since $j(\mu) > j'' \mu$, there
is an ordinal $\gb < j(\mu)$ above all the ordinals in $j  ''\mu$.
 Now, for each $\ga < \mu$ we can find some $\gz< \chi$ such that
  $$j(F)(j(\ga) , \gb)=\gz.$$
  Since  $\chi < \gl_0$, and as
 $\mu$ is regular in $V$ and no new sequences of length $< \lambda_0$ are
added to $V$ in $V[G]^P$,
we may find a single
$\gz$ such that for unboundedly many $\ga$'s the
equality $j(F)(j(\ga),\gb)=\gz$ holds.
Since $j$ is elementary, it follows that for any $\ga_1 < \ga_2$ in this
 unbounded set
$F(\ga_1,b)=F(\ga_2,b)=\gz$ holds for some $b > \ga_2$.  Thus
an unbounded  $\mu$--branch for $\gz$ was found in $V[G]^P$,
which is a $\gl_0$--closed extension of $V[G]$.

We turn now to the proof of Lemma \ref{extension}.
The collapsing poset
 $C_1= \Coll(\lambda_0, <\lambda_1)$ is a factor of
$C=\Coll(\seqlambda)$, and for simplicity of expression, we identify $c \in
C_1$ with the condition $\la c_1 \mid i < \omega \ra \in C$ defined by
$c_0=c$ and $c_i= \emptyset$ for $i>0$.

Denote each $j(\gl_i)$ with $\gl_i^*$.  Then
$\gl^*_0=\gl_0$, but $\gl_1^* > \mu$.  In $M$, $j(C)$ is $[\Coll(\la
\gl_i^* \mid i \in \go \ra )]^M$, and
 $C_1^*=\Coll(\gl_0,< \gl_1^*)$ (which is the same---defined in $V$
or in $M$), is a factor of $j(C)= [\Coll(\la \gl_i^* \mid i< \omega
 \ra )]^M$.

Let $G \subseteq C$ be a $V$-generic filter over $C$. Observe that if $\la
c_i \mid i < \omega \ra \in G$, then $c \in G$ as well. In order to extend
$j$ on $V[G]$ and to prove the lemma, we should find in $V[G]$ a
$\lambda_0$-complete poset $P$ such that in $V[G]^P$ there is a $V$-generic
filter $G^\ast$ over $j(C)$ such that
$$\mbox{If}\ g \in G,\ \mbox{then}\ j(G) \in G^\ast.$$
If we do so, then an embedding of
$V[G]$ into $M[G^\ast]$ can be defined as follows: For any $x \in V[G]$, let
${\bf x}$ be a name of $x$ in $V^C$. Then $j({\bf x})$ is a name in
$V^{j(C)}$ and we define $j'(x)$ to be its interpretation in $V[H']$.
We trust the reader to check that $j'$ is a well defined elementary
extension of $j$.

nstead of writing down $P$, we will describe an iteration of two
extensions, each one $\lambda_0$-complete.

Since
 $\mu < \gl_1^*$ is collapsed to $\lambda_0$ in $V^{C_1^*}$, Lemma
\ref{l22}
implies that there is a projection, $\Pi$, of $C_1^*$ onto $C$, which
can be used to find a generic extension of $V[G]$ which has the form $V[H]$
for a $V$-generic filter $H$ over $C_1^\ast$ such that:
\begin{enumerate}
\item The passage from $V[G]$ to $V[H]$ is done by forcing with a
$\lambda_0$-closed forcing.
\item For every $c \in C_1$, $c \in G$ iff $c \in H$.
\end{enumerate}
Thus, for every $g \in G$, $j(g) \in j(C)/H$. Indeed, any $g \in G$ has the
form $g= \la c, \overline{r} \ra$ where $c \in C_1 \cap G$, and
$\overline{r}$ is the later part of the sequence. Then $j(g) = \la c,
j(\overline{r}) \ra$ where $c \in C_1^\ast \cap H$, and thus $j(g) \in
j(C)/H$.

It follows that $j''G$ (the image of $G$ under the restriction of $j$ to
$C$) is in $M[H]$ a pairwise compatible collection of conditions in
$j(C)/H$. Since $j(C)/H$ is isomorphic to $\Coll(\la \lambda^\ast_i \mid i
\geq 1 \ra)$ in $M[H]$, it is $\lambda_1^\ast$-complete, and a supremum,
denoted $s \in j(C)/H$ can be found for $j''G$. This is our ``master
condition'': If $G^\ast$ is any $V[H]$-generic filter over $j(C)/H$,
containing $s$, then:
\begin{enumerate}
\item $G^\ast$ is in fact $V$-generic over $j(C)$.
\item The forcing $j(C)/H$ is $\lambda_1^\ast$-complete in $M[H]$, and it
is thence $\lambda_0$-complete in $V[H]$ (because $M[H]$ is
$\lambda_0$-closed in $V[H]$).
\end{enumerate}

\section{Conclusion}
We have proved the following theorem:
\begin{theorem}
Assume a cardinal $\kappa$ and sequence $L=\seqlambda$ such that
\begin{description}
\item[$B_1$:] For $\lambda=\sup \{ \lambda_i \mid i < \omega \}$ and $\mu =
\lambda^+$, each $\lambda_i$, $ i > 0$, is $\mu$-supercompact.
\item[$B_2$:]$\kappa$ is the critical point of an embedding $j : V
\longrightarrow M$ where $j(\kappa) = \lambda_0$ and $M^\mu \subseteq M$.
\end{description}
Then there is a generic extension in which there are no $\alephomp$ \A\
trees.
\end{theorem}
Indeed, in Section \ref{s22} we saw that by making a preparatory extension
we may assume that $\kappa$ is such that if $C= \Coll(L)$, then in $V^C$
$\kappa$ is $\kpompo$-huge. So, we go to $V^C$, and find that the
potential branching property for $\mu=\jkpompo$ holds (by Theorem
\ref{t61}). But now, in $V^C$, all the assumptions for theorem \ref{notrees}
hold. Thus, in a final extension, obtained as a product of $\Coll(\kpompo,
< j(\kappa))$ and $\Coll(\aleph_0, \kpom)$, there are no \A\ trees on
$\alephomp$.
\section*{References}
Baumgartner, Iterated forcing (on reverse Easton)?\\
Jensen, Weak square = special tree\\
Magidor, On reflecting stationary sets\\
Magidor Levinski Shelah, Chang's conjecture on $\alephom$\\
Schimmerling no weak square implies inner model with Woodin.\\
Solovay Strong axioms, GCH above supercompacts.
\end{document}